\documentclass{article}

\usepackage{lineno,hyperref}
\usepackage{epsfig}
\usepackage{graphicx}
\usepackage{epstopdf}
\usepackage{amssymb,amsmath,amsthm}
\usepackage{color}

\newtheorem{example}{Example}

\newtheorem{proposition}{Proposition}
\newtheorem{corollary}{Corollary}

\begin{document}

\title{Solving differential and integral equations with Tau method\footnote{This work was partially supported by CMUP(UID/MAT/00144/2013), which is funded by FCT (Portugal) with national and European structural funds (FEDER), under the partnership agreement PT2020}}

\author{Jo\~ao Carrilho de Matos\footnote{jem@isep.ipp.pt, Instituto Superior de Engenharia do Porto, Rua Dr. Ant\'onio Bernardino de Almeida, 431, 4249-015 Porto, Portugal}, Jos\'e M. A. Matos\footnote{Instituto Superior de Engenharia do Porto and Centro de Matem\'atica da Universidade do Porto}, Maria Jo\~ao Rodrigues\footnote{Faculdade de Ci\^encias da Universidade do Porto and  Centro de Matem\'atica da Universidade do Porto}}

\maketitle

\begin{abstract}
In this work we present a new approach for the implementation of operational Tau method for the solutions of linear differential and integral equations. In our approach we use the three terms relation of an orthogonal polynomial basis to compute the operational matrices. We also give numerical applications of operational matrices to solve differential and integral problems using the operational Tau method.
\end{abstract}

\textbf{keywords:} Operational Tau method; Orthogonal polynomials; Differential equations; Integral equations

\section{Introduction}
\label{intro}
The operational Tau method, \cite{Os80}, is a spectral method for solving differential and integral equations. 
These methods use matrices (called operational matrices) to represent linear operators defined
in linear function spaces in a given orthogonal basis (see for instance \cite{Can87,Fu92,Go77}).
The original operational Tau approach has a serious drawback since the operational matrices are computed using
similar matrices, 
which are ill-conditioned. Thus, this approach is not suitable for problems that require high order approximations.

In this paper we avoid the use of similar matrices 
 building the operational matrices  using
the tree therm recurrence relation associated to a given orthogonal polynomial basis. 
We also give numerical examples applying our approach to integral and differential equations.

\section{Operational Tau method}
\label{sec:1}
The key idea of the operational Tau method formulation, given in \cite{Os80} and \cite{Os81}, is to represent in matrix form linear differential operators with polynomial coefficients. This matrix representation can be generalized to integral or integro-differential operators. 

\subsection{Matrix representation of linear operators in power basis}
Let $\mathbb{P} [x]$ and $\mathbb{P}_{n} [x]$ denote the linear space of polynomials and the linear space of polynomials of degree at most $n$
in one variable, $x$, respectively.  
Let $\mathcal{L}$ be a linear differential operator with polynomial coefficients
\begin{equation}\label{linopdef}
\mathcal{L}=\sum_{i=0}^{m}p_{i}(x)\frac{\mathrm{d}^{i}}{\mathrm{d}x^{i}}, \ \ \ p_{i}(x)=\sum_{j=0}^{n_{i}}p_{i,j}x^{i},\in\mathbb{P}_{n_{i}} [x],
\end{equation}
\noindent and let $y_{n}(x)\in \mathbb{P}_{n} [x]$, $y_{n}(x)=\sum_{i=0}^{n}a_{i}x^{i}$ written as
$\mathbf{ y}_{n}=\mathbf{ a}\, \mathbf{x}$ in a matrix form. Then $\mathcal{L}\left[y_{n}(x)\right]$ has the following matrix representation, in the power basis,
$$
\mathcal{L}\left[y_{n}(x)\right]=\mathbf{ a}\,\boldsymbol{\Pi} \, \mathbf{x},
$$
\noindent where the matrix $\boldsymbol{\Pi}$ is defined by $\displaystyle \boldsymbol{\Pi}=\sum_{i=0}^{m} p_{i}\left( \mathbf{ M} \right)\mathbf{ H^{i}}$, with matrices $\mathbf{ H}$ and $\mathbf{M}$ representing
the linear differential and shift operator, respectively.  That is, we have
\begin{eqnarray*}
\frac{\mathrm{d}^{k}}{\mathrm{d}x^{k}}\left[ y_{n}(x)\right]= \mathbf{ a}\, \mathbf{ H}^{k}\, \mathbf{ x} \\
x^{k}\,  y_{n}(x)= \mathbf{ a}\, \mathbf{M}^{k}\, \mathbf{ x} 
\end{eqnarray*}
\noindent with
$$
\mathbf{ H}=
\left[
\begin{array}{ccc}
     0 &   &   \\
     1 & 0&   \\
     &2 &0 \\
     &  & \cdots
\end{array}
\right],
\ \ \ 
\mathbf{M}=
\left[
\begin{array}{cccc}
     0 &  1 &  & \\
     & 0& 1 & \\
     &&0 &1\\
     &  & &\cdots
\end{array}
\right].
$$

\noindent We may also  generalize this matrix representation of the operator  $(\ref{linopdef})$ to integral 
linear operators, with polynomials coefficients, using the matrix
$$
\boldsymbol{\Theta}=
\left[
\begin{array}{cccc}
     0 &  1 &  & \\
     & 0& \frac{1}{2} & \\
     &&0 &\frac{1}{3}\\
     &  & &\cdots
\end{array}
\right].
$$
In fact, the primitive of polynomial $y_{n}$ (with zero constant) can be given by
$$
\int y_{n}(x)\, \mathrm{d}x = \mathbf{ a}\,\boldsymbol{\Theta}\,\mathbf{ x}.
$$


\subsection{Classic approach of operational Tau method}

Now, let us consider a matrix $\boldsymbol{\nu}=[\nu_{0}, \nu_{1}, \nu_{2},\cdots]^{T}$, where
$\nu=(\nu_{0}, \nu_{1}, \nu_{2},\cdots)$ is a basis of $\mathbb{P}[x]$ such that, for each non negative integer $k$, $\nu_{k}$ is a polynomial of degree $k$. Thus, the image of the polynomial $y_{n}$, expanded on basis $\nu$, of the operator $\mathcal{L}$ is given by 
$$
\mathcal{L}[y_{n}]=\mathbf{ a}\, \boldsymbol{\Pi_{\nu}}\, \boldsymbol{\nu},
$$
\noindent where, $\boldsymbol{\Pi_{\nu}}=\mathbf{V}\, \boldsymbol{\Pi}\, \mathbf{ V}^{-1}$, and 
$\mathbf{ V}$ is the change of basis matrix that satisfies the relation $\boldsymbol{\nu}=\mathbf{ V}\,\mathbf{ x}$.

Consider a linear problem
\begin{eqnarray}\label{problem1}
\mathcal{L}u& = & f, \ \ \ x\in]a,b[ \\
g_{j}(u) & = & \sigma_{j}, \ \ \ j=1,2,\cdots, m 
\end{eqnarray}
\noindent where, $g_{j}$,  $ j=1,2,\cdots, m$ are linear functionals that represent the supplementary conditions and
$f\in\mathbb{P}[x]$. A Tau solution of order $n$ expanded on a basis $\nu$ is the solution $u_{n}\in\mathbb{P}_{n}[x]$ of the associated problem to
the problem (\ref{problem1}) 
\begin{eqnarray}\label{problem1}
\mathcal{L}u& = & f+h_{n}, \ \ \ x\in]a,b[ \\
g_{j}(u) & = & \sigma_{j}, \ \ \ j=1,2,\cdots, m 
\end{eqnarray}
\noindent where $h_{n}$ it is a perturbation polynomial. The coefficients, $a_{i}$, $i=0,1,\cdots, n$ of the tau solution $u_{n}(x)=\sum_{i=0}^{n}a_{i}\nu_{i}(x)$ are solutions of a system of a linear equations \cite{Os80}. 

We note that this system includes the matrix, $\boldsymbol{\Pi_{\nu}}$, that represent the linear operator $\mathcal{L}$, the  supplementary conditions and the polynomial $f$. Ano\-ther important remark is that we use infinity matrices to represent the linear operators. However, we need not to worry about the meaning of 
matrix multiplication. From the practical point of view we deal with polynomials. Thus all this products reduce to a finite number of non null parcels and the size of the finite matrices that we work depend on the number of supplementary conditions and on the hight of operator $\mathcal{L}$. For details see \cite{Os80} and \cite{Os81}.

\section{Matrix representation of linear operators in orthogonal basis}

Consider an orthogonal polynomial basis $\nu=(\nu_0,\nu_1,\ldots)$ of $\mathbb{P}[x]$ defined by an inner product
\begin{equation}\label{PiPj}
<\nu_i,\nu_j>=\int^{b}_{a} \nu_i\ \nu_j\ w \mathrm{d}x = ||\nu_i||^2\delta_{i,j},\ i,j\in\mathbb{N}_0 ,
\end{equation}
\noindent where, as usual,  $w$ is a weight function and $||*||$ is the associated norm to the inner product $<*,*>$.

The Fourier coefficients $f_{i}$ of a function $f$ are given by
  \begin{equation}\label{fi}
	 f_i=\frac{1}{||\nu_i||^2}<\nu_i,f> 
 \end{equation}
\noindent Assuming  that all integrals $\int^{b}_{a} \nu_i\ \nu_j\ w \mathrm{d}x$ exist we will write, $f(x)=\sum_{i=0}^{\infty}f_{i}\nu_i(x)$, where the equality only holds when the infinite series converge to $f$. Then we have the following

\begin{proposition}\label{prop:L}
	If $\mathcal{L}$ is a linear operator acting on $\mathbb{P}$ and $\mathbf{ L}_{\nu}$ is the infinite matrix defined by
	\begin{equation}
	\mathbf{ L}_{\nu}=[\ell_{i,j}]_{i,j\geq 0},\quad \text{with}\quad \ell_{i,j} = \frac{1}{||\nu_i||^2} <\nu_i,\mathcal{L}[\nu_j]>,
	\end{equation}
	then formally $\mathcal{L}{\nu}={\boldsymbol{\nu}}\, \mathbf{ L}_{\nu}$.
	
	\begin{proof}
		For each $j=0,1,\ldots$ we define the infinite unitary vector $\mathbf{ e}_j=[\delta_{i,j}]_{i\geq 0}$. So that ${\nu} e_j = \nu_j$ and using \eqref{fi} we get
		\[ \mathcal{L}[\boldsymbol{\nu} \mathbf{ e}_j ]= \mathcal{L}[\nu_j ]= \sum_{i\geq 0} \frac{1}{||\nu_i||^2} <\nu_i, \mathcal{L}[\nu_j]> \nu_i = \sum_{i\geq 0} \ell_{i,j} \nu_i = \boldsymbol{\nu} \mathbf{ L}_{\nu} \mathbf{ e}_j,\ j=0,1,\ldots  \]
		and so  $\mathcal{L}{\nu}={\boldsymbol{\nu}}\, \mathbf{ L}_{\nu}$, in the element wise sense.
	\end{proof} 
\end{proposition}

It is well known that a sequence of orthogonal polynomials, normalized with the condition $\nu_{0}=1$, satisfies a three term recurrence relation.

\begin{equation}\label{ttrr}
\left\{\begin{array}{l}
x \nu_j = \alpha_j \nu_{j+1} + \beta_j \nu_{j} + \gamma_j \nu_{j-1},\ j\geq 0 \\
\nu_{-1}=0,\ \nu_0=1
\end{array}\right. ,
\end{equation} 

\noindent This recurrence relation it is useful to  find the matrices $\mathbf{ M}_{\nu}$, $\mathbf{H}_{\nu}$ and $\boldsymbol{\Theta}_{\nu}$ that represent, respectively, the shift, differential and integral operators in $\nu$ basis \cite{JCM17}.

For the shift operator we have,
\begin{proposition}\label{prop:xP}
	Let $\nu$ be a basis satisfying \eqref{ttrr}, defining 
	\[
	\mathbf{ M}_{\nu}=[\mu_{i,j}]_{i,j\geq 0},\quad \text{with}\quad \mu_{i,j} = \frac{1}{||\nu_i||^2}<\nu_i,x \nu_j>,
	\]
	then 
	\begin{equation}\label{mu} 
		\begin{cases}
		\mu_{0,0} = \beta_0,\ \mu_{1,0} =  \alpha_0 \\
		\mu_{j-1,j} = \gamma_j,\  \mu_{j,j} = \beta_j,\ \mu_{j+1,j} =  \alpha_j \\ 
		\mu_{i,j} = 0,\ |i-j|>1 
		\end{cases},\ j=1,2,\ldots	
	\end{equation}
	and  $x {\nu}={\boldsymbol{\nu}} \mathbf{ M}_{\nu}$.
	\begin{proof}
		From the definition of $\mu_{i,j}$ and \eqref{fi} follows that
		\[ x\nu_j=\sum_{i=0}^{j+1} \mu_{i,j} \nu_i,\ j=0,1,\ldots \]
		and using \eqref{ttrr} we get \eqref{mu}. The fact that $x {\nu}={\boldsymbol{\nu}} \mathbf{ M}_{\nu}$, in the element-wise sense, is a consequence of proposition \ref{prop:L}.
	\end{proof} 
\end{proposition}

For the differential operator stands the following

\begin{proposition}\label{prop:dP}
	Let $\nu$ be a basis satisfying \eqref{ttrr}, defining 
	\[
	\mathbf{ H}_{\nu} = [\eta_{i,j}]_{i,j\geq 0},\quad \text{with}\quad \eta_{i,j} = \frac{1}{||\nu_i||^2}< \nu_i, \frac{\mathrm{d}}{\mathrm{d}x} \nu_j >,
	\]
	then  for $j=1,2,\ldots$
	\begin{equation}\label{eta} 
	 	\begin{cases}
		\eta_{i,j+1}=\frac{1}{\alpha_j}\big[\alpha_{i-1}\eta_{i-1,j} + (\beta_i-\beta_j)\eta_{i,j} \\ \hspace{0.24\textwidth} + \gamma_{i+1}\eta_{i+1,j} - \gamma_j\eta_{i,j-1}\big],\ i=0,\ldots,j-1\\
		\eta_{j,j+1}=\frac{1}{\alpha_j}(\alpha_{j-1}\eta_{j-1,j}+1)
		\end{cases},
	\end{equation}
	and  $\frac{\mathrm{d}}{dx} {\nu}={\boldsymbol{\nu}} \mathbf{ H}_{\nu}$.
	\begin{proof}
		Applying the operator $\frac{\mathrm{d}}{\mathrm{d}x}$ to both sides of \eqref{ttrr} then
		\[ \nu_j + x \nu'_j = \alpha_j \nu'_{j+1} + \beta_j \nu'_{j} + \gamma_j \nu'_{j-1},\ j=0,1,\ldots  \]
		and, by definition of $\eta_{i,j}$
		\[ \nu_j + x \sum_{i=0}^{j-1}\eta_{i,j}\nu_i = \alpha_j \sum_{i=0}^{j}\eta_{i,j+1}P_i + \beta_j \sum_{i=0}^{j-1}\eta_{i,j}\nu_i + \gamma_j \sum_{i=0}^{j-2}\eta_{i,j-1}\nu_i,\ j=0,1,\ldots  \]
		and so
		\begin{eqnarray*}
		\alpha_j \sum_{i=0}^{j}\eta_{i,j+1}\nu_i  & = &  \nu_j + \sum_{i=0}^{j-1}\eta_{i,j} (\alpha_i \nu_{i+1} + \beta_i \nu_{i} + \gamma_i \nu_{i-1}) \\
		& & \qquad\qquad\qquad - \beta_j \sum_{i=0}^{j-1}\eta_{i,j} \nu_i - \gamma_j \sum_{i=0}^{j-2}\eta_{i,j-1}\nu_i 
		\end{eqnarray*}
		rearranging indices and identifying similar coefficients we get \eqref{eta}. That $\frac{\mathrm{d}}{\mathrm{d}x} {\nu}={\boldsymbol{\nu}} \mathbf{ H}_{\nu}$ is a consequence of Proposition \ref{prop:L}.
	\end{proof} 
\end{proposition} 

To derive the matrix $\boldsymbol{\Theta}_{\nu}$ we have,

\begin{proposition}\label{prop:IP}
	Let $\nu$ be the basis satisfying \eqref{ttrr}, defining  
	\[
	 \boldsymbol{\Theta}_{\nu} = [\theta_{i,j}]_{i,j\geq 0},\quad \text{with}\quad \theta_{i,j} = \frac{1}{||\nu_i||^2}<\nu_i, \int \nu_j \, \mathrm{d}x >,
	\]
	then  for $j=1,2,\ldots$
	\begin{equation}\label{thetaij}
	\left\{\begin{array}{l}
	\theta_{j+1,j} = {\displaystyle \frac{\alpha_j}{j+1}} \\
	\theta_{i+1,j} = {\displaystyle -\frac{\alpha_i}{i+1} \sum_{k=i+2}^{j+1} \eta_{i,k}\theta_{k,j}}, i=j-1,\ldots,1,0
	\end{array}\right. 
	\end{equation}
	and  $\int {\nu}\mathrm{d}x={\boldsymbol{\nu}} \boldsymbol{\Theta}$.
	\begin{proof}
		By definition, considering that the primitive of $\nu_j$ is a polynomial of degree $j+1$ defined with an arbitrary constant term, we can write
		\[ \int \nu_j \mathrm{d}x = \sum_{i=1}^{j+1}\theta_{i,j}P_i,\ j=0,1,\ldots  \]
		Differentiating both sides and applying proposition \ref{prop:dP} 
		\[
		\nu_j = \sum_{i=1}^{j+1}\theta_{i,j}\nu'_i = \sum_{i=1}^{j+1}\theta_{i,j}\sum_{k=0}^{i-1}\eta_{k,i} \nu_k .
		\]
		Rearranging indices and identifying similar coefficients, 
		\[
		\nu_j = \sum_{i=0}^{j}\left[ \sum_{k=i+1}^{j+1} \eta_{i,k} \theta_{k,j}  \right] \nu_i .
		\]
		And so, for the coefficient of $\nu_j$,
		\[ 
		\eta_{j,j+1}\theta_{j+1,j} = 1
		\]
		and, for the coefficients of $\nu_i,\ i=0,\ldots,j-1$,
		\[
		\sum_{k=i+1}^{j+1}\eta_{i,k}\theta_{k,j} = 0
		\]
		The result is obtained solving for $\theta_{j+1,j}$ the first equation and for $\theta_{i+1,j}$ each one in the last set of equations.
		
	\end{proof} 
\end{proposition}

\begin{corollary}\label{corol:IPvolterra}
Let $\nu$ be the basis satisfying \eqref{ttrr}, $y=\boldsymbol{\nu} \mathbf{ a}$ a formal Fourier series and the matrix $\boldsymbol{\Theta}_{\nu}$ defined in
Proposition \eqref{prop:IP}. Then $\int_{a}^{x}y=\boldsymbol{\nu}\boldsymbol{\Theta}_{\nu}^{x,a}\mathbf{ a}$ where,
$$
\boldsymbol{\Theta}_{\nu}^{x,a}=[v_{i,j}]_{i,j\geq0}, \ \
\left\{
\begin{array}{l}
v_{0,j}=-\sum_{i=1}^{j+1}\theta_{i,j}\nu_{i}(a)\\
\\
v_{i,j}=\theta_{i,j}, \ i>0
\end{array} 
\right.
$$

\begin{proof} 
Defining the primitive with undefined $\nu_{0}$, $F_{j}=\int \nu_{j}=\sum_{i=1}^{j+1}\theta_{ij}\nu_i$, then
		\begin{eqnarray*}
			 \int_{a}^{x} \nu_j(t)dt & = & F_j(x) - F_j(a) \\
			 & = &  \sum_{i=1}^{j+1}\theta_{i,j}\nu_i(x) - \sum_{i=1}^{j+1}\theta_{i,j}\nu_i(a) = \nu \boldsymbol{\Theta}_{\nu}^{x,a} e_j
		\end{eqnarray*}
		meaning that $\int_{a}^{x}\boldsymbol{\nu} = \boldsymbol{\nu} \boldsymbol{\Theta}_{\nu}^{x,a} $ in element wise sense. The proof follows by linearity. 

\end{proof} 
\end{corollary}

\section{Numerical results}

In order to test the numerical stability of the recurrence relations, used to compute the matrices, $\boldsymbol{M}_{\nu}$ and $\boldsymbol{\Theta}_{\nu}$, we will solve three stiff problems using the operational Tau method. We have chosen two differential problems and a integral 
one to test the recurrence relations in Jacobi  and Laguerre polynomials cases. In all problems the  Tau method converges slowly and we need to compute the higher order operational matrices to obtain a good approximation of the solution of the given problem. Thus, if we can improve the accuracy of higher order Tau solutions, we can conclude (since  the matrix $\boldsymbol{\Pi}$
is not ill conditioned) that the recurrence relations used for the orthogonal polynomials, are numerically stable.

In the next two examples it will be useful the three terms relation \eqref{ttrr}, see e.g. \cite{AbSteg}, for
Jacobi polynomials $P^{\alpha,\beta}_{n}(x)$, $n\geq0$ with \newline
	$\alpha_n= \frac{2(n+1)(n+\gamma+1)}{(2n+\gamma+1)(2n+\gamma+2)},\
	\beta_n = \frac{(\beta-\alpha)\gamma}{(2n+\gamma)(2n+\gamma+2)},\
	\gamma_n = \frac{2(n+\alpha)(n+\beta)}{(2n+\gamma)(2n+\gamma+1)}$ \newline
where $\alpha,\ \beta>-1$ and $\gamma=\alpha+\beta$.
 
\begin{example}
We consider the following differential problem with boundary conditions
 
\begin{equation}\label{eq:turningpoint}
\left\{\begin{array}{l}
\epsilon y''(x)-xy(x)=0,\ x\in ]-1,1[ \\
y(-1) = 1, \ y(1)=1
\end{array}\right. 
\end{equation}
where $\epsilon$ is a real positive parameter.
This problem has solution $$y(x)=c_{1}\mathrm{Ai}\left(\frac{x}{\sqrt[3]{\epsilon}}\right)+c_{2}\mathrm{Bi}\left(\frac{x}{\sqrt[3]{\epsilon}}\right),$$ where $\mathrm{Ai}$ and $\mathrm{Bi}$ are the airy functions of first and second kind, respectively, and the constants $c_{1}$ and $c_{2}$ are computed in such a way to fulfill the boundary conditions.
For small values of the parameter $\epsilon$ the solution has a smooth and a strong oscillations regions, see Figure \ref{fig:turningpoints1}.

We solve this problem with $\epsilon=10^{-5}$ using Jacobi polynomials $P^{(\alpha,\beta)}_{n}(x)$ as basis.
For this example the operational matrix   is given by
$$
\boldsymbol{\Pi}_{P^{(\alpha,\beta)}}=\epsilon \boldsymbol{H}_{P^{(\alpha,\beta)}}^{2}-\boldsymbol{M}_{P^{(\alpha,\beta)}}.
$$ 

\begin{figure}
	\includegraphics[width=\textwidth]{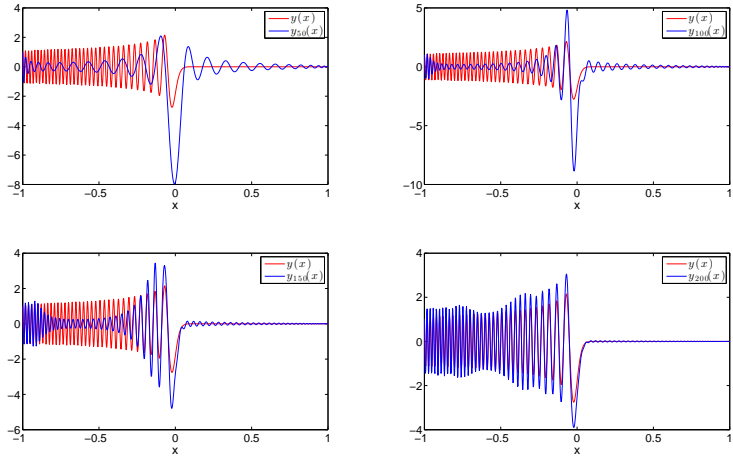}
	\caption{Blue line: Jacobi-Tau approximations $y_{50}$, $y_{100}$, $y_{150}$ and $y_{200}$, with $\alpha=\beta=0$ (Legendre basis). Red line: solution of the problem \eqref{eq:turningpoint}}
	\label{fig:turningpoints1}
\end{figure}

\begin{table}[hbp]
\centering\large
\caption{Infinity norm of Jacobi-Tau approximation errors, $max_{x\in[-1,1]}|y(x)-y_{n}(x)|$ of equation (\ref{eq:turningpoint}).}	
\begin{tabular}{ccccc}
&\multicolumn{4}{ c }{$\mathbf{ n}$} \\
\hline
$\boldmath{(\alpha, \beta)}$ &$\mathbf{ 150}$&$\mathbf{ 250}$&$\mathbf{ 350}$&$\mathbf{ 1000}$\\ \hline \hline
$(0,0)$                                &  $2.05e+0$  &  $3.32e-8$ & $3.48e-12$ & $ 5.29e-12$\\ 
$(-1/2,-1/2)$                       &  $6.42e+0$  &  $4.11e-8$ & $6.04e-12$ & $ 5.77e-12$\\ 
$(1,-9/10)$                             &  $6.46e+1$  &  $3.41e-8$ & $2.44e-11$ & $ 2.86e-11$\\ 
$(-9/10,-9/10)$                &  $2.37e+1$  &  $4.88e-8$ & $1.86e-11$ & $ 3.96e-12$\\ 
$(1/2,-1/2)$                    &  $1.75e+0$  &  $2.68e-8$ & $2.12e-12$ & $ 1.17e-12$\\ 
\end{tabular}
\label{tab:prob1}
\end{table}

We can see on Table \ref{tab:prob1} that this Tau problem has slow rate of convergence. In fact we need a Tau approximation of degree $250$ to reach an error of order $10^{-8}$ and an approximation of degree $350$ to reach errors of order $10^{-11}$ or $10^{-12}$ depending on the Jacobi basis (i. e. depending on the values of $\alpha$ and $\beta$). Thus, the computed recurrence relations, for the derivative, given in Proposition \ref{prop:dP} , are stable.
The solutions of higher orders are also good approximations. In fact, for 
$n>350$ the matrices $\boldsymbol{\Pi}_{P^{(\alpha,\beta)}}$ are ill conditioned implying that it is useless to increase the degree of the Jacobi -Tau approximation. 
\end{example}


In the following example we  test the stability of integral recurrence relation given on proposition \ref{prop:IP}

\begin{example}

Consider the Volterra integral equation

\begin{equation}\label{volterraeq}
(x-a)^{3}y(x)+\int_{-1}^{x}y(s)ds=-f(-1), \ \ \ x\in[-1,1]
\end{equation}
\noindent where $f(x)=\exp\left(\frac{1}{2(x-a)^{2}}\right)$ and $a$ is a real parameter. The solution of \eqref{volterraeq} it is 
the function $y$ defined by
$
y(x)=(a-x)^{-3}f(x).
$

We solve this problem with $a=1.25$ using Jacobi polynomials $P^{(\alpha,\beta)}_{n}(x)$ as basis.
For this case the operational matrix is given by
$$
\boldsymbol{\Pi}_{P^{(\alpha,\beta)}}= (\boldsymbol{M}_{P^{(\alpha,\beta)}}-a \boldsymbol{I})^{3}+\boldsymbol{\Theta}_{P^{(\alpha,\beta)}}^{x,-1},
$$
\noindent where the matrix $\boldsymbol{\Theta}_{P^{(\alpha,\beta)}}^{x,-1}$ is given on Corolary \ref{corol:IPvolterra}.

\begin{table}[hbp]
\centering\large
\caption{Infinity norm of Jacobi-Tau approximation errors, $max_{x\in[-1,1]}|y(x)-y_{n}(x)|$ of equation (\ref{volterraeq}).}	
\begin{tabular}{ccccc}
&\multicolumn{4}{ c }{$\mathbf{ n}$} \\
\hline
$\boldmath{(\alpha, \beta)}$ &$\mathbf{ 50}$&$\mathbf{ 100}$&$\mathbf{ 150}$&$\mathbf{ 1000}$\\ \hline \hline
$(0,0)$                                &  $3.90e+0$  &  $1.85e-7$ & $1.58e-7$ & $ 1.59e-7$\\ 
$(-1/2,-1/2)$                       &  $1.30e+0$  &  $5.42e-7$ & $5.46e-7$ & $ 5.46e-7$\\ 
$(1,-9/10)$                             &  $3.49e+1$  &  $5.26e-7$ & $4.02e-9$ & $ 3.78e-9$\\ 
$(10,0)$                &  $1.57e+4$  &  $7.35e-2$ & $1.34e-9$ & $ 1.72e-9$\\ 
\end{tabular}
\label{tab:prob2}
\end{table}

The results presented on Table \ref{tab:prob2} show that the behavior of the Tau solutions of this example is similar to the previous one. We need higher order Tau solutions to reach errors of order $10^{-7}$ or $10^{-9}$ (depending on the Jacobi polynomial basis) and our approach  is stable. 
\end{example}

In the next example we analyze the behavior of the recurrence relation given on Proposition \ref{prop:dP} for Laguerre
polynomials basis. 
The Laguerre polynomials $L_{n}(x)$, $n\geq0$ satisfy \eqref{ttrr} with
\begin{equation*} 
\alpha_n=-(n+1),\ \beta_n=2n+1,\ \gamma_n = -n.
\end{equation*} 
\begin{example}
Here we consider the Bessel equation with boundary conditions
\begin{equation}\label{eq:bessel}
\left\{\begin{array}{l}
 x^2 y''(x)+xy'(x)+(x^2-m^2)y(x)=0,\ x\in ]0,60[ \\
y(0) = 0, \ y(60)=1
\end{array}\right. 
\end{equation}
\noindent where $m\neq 0$ it is a real parameter. This problem has solution
$$
y(x)=\frac{\mathrm{J}_{m}(x)}{\mathrm{J}_{m}(60)}
$$
\noindent where $\mathrm{J}_{m}(x)$ is the Bessel function of first kind.

\begin{figure}
	\includegraphics[width=\textwidth]{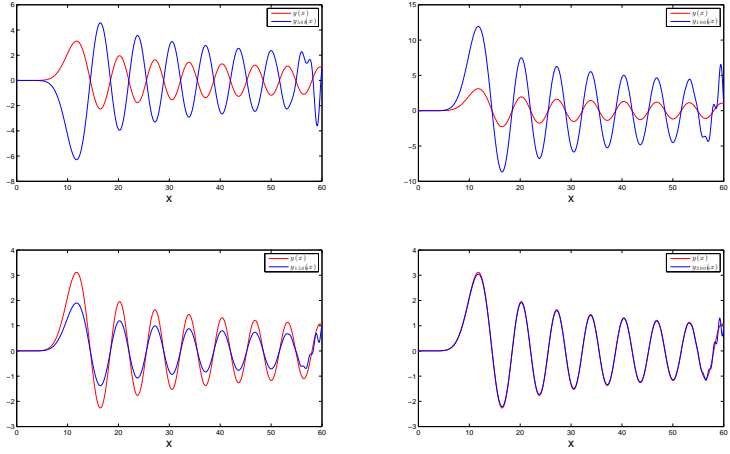}
	\caption{Blue line: Laguerre -Tau approximations $y_{500}$, $y_{1000}$, $y_{1500}$ and $y_{2000}$, with $m=10$. Red line: solution of the problem \eqref{eq:bessel}.}
	\label{fig:laguerre1}
\end{figure}
\end{example}

 We show on Figure \ref{fig:laguerre1} the graphs of the Laguerre-Tau approximations (blue lines), $y_{500}$, $y_{1000}$, $y_{1500}$ and $y_{2000}$ of this problem with the parameter value $m=10$. We can see that although the rate of convergence   is very slow, it is possible reach good approximations using our approach for the operational Tau method. It is still possible to compute approximations with degree higher than $2000$ but that does not improve significantly
the accuracy obtained by $y_{2000}$.
  
\section{Conclusions}
Numerical results show that the recurrence process to build operational matrices applied to operational Tau method stabilizes the ''classic'' Tau method, introduced in \cite{Os80}. More, our approach allows to work in all orthogonal polynomial bases.

\end{document}